\documentclass[a4paper, 12pt]{article}

\usepackage{graphicx}
\usepackage{url}
\usepackage{amsbsy,amssymb} 
\usepackage{amsmath}
\usepackage{abstract}
\usepackage{hyperref}
\usepackage{color}

\usepackage{amsmath,amsthm}

\usepackage[all]{xy} 

\newtheorem{Proposition}{Proposition}
\newtheorem{Corollary}{Corollary}

\newtheorem{Theorem}{Theorem}

\newtheorem{Lemma}{Lemma}
\newtheorem{Definition}{Definition}

\begin{document}

{\LARGE\centering{\bf{The Poincare lemma, antiexact forms, and fermionic quantum harmonic oscillator}}}

\begin{center}
\sf{Rados\l aw Antoni Kycia$^{1,2,a}$}
\end{center}

\medskip
\small{
\centerline{$^{1}$Masaryk University}
\centerline{Department of Mathematics and Statistics}
\centerline{Kotl\'{a}\v{r}sk\'{a} 267/2, 611 37 Brno, The Czech Republic}
\centerline{\\}
\centerline{$^{2}$Cracow University of Technology}
\centerline{Faculty of Materials Engineering and Physics}
\centerline{Warszawska 24, Krak\'ow, 31-155, Poland}
\centerline{\\}

\centerline{$^{a}${\tt
kycia.radoslaw@gmail.com} }
}

\begin{abstract}
\noindent
The paper focuses on various properties and applications of the homotopy operator, which occurs in the Poincar\'{e} lemma. In the first part, an abstract operator calculus is constructed, where the exterior derivative is an abstract derivative and the homotopy operator plays the role of an abstract integral. This operator calculus can be used to formulate abstract differential equations. An example of the eigenvalue problem that resembles the fermionic quantum harmonic oscillator is presented. The second part presents the dual complex to the Dolbeault bicomplex generated by the homotopy operator on complex manifolds.
\end{abstract}
Keywords: Poincare lemma, antiexact differential forms, homotopy operator, fermionic harmonic oscillator, complex manifold; \\
Mathematical Subject Classification: 58A12, 58Z05;\\

\section{Introduction}

The Poincar\'{e} lemma is one of the most important tools of exterior calculus. Although it is a very old result it is continuously generalized in various ways \cite{PoincareGeneralization1, PoincareGeneralization2, PoincareGeneralization3, DiscretePoincareLemma, DiscreteExteriorCalculus}, including non-Abelian cases \cite{PoincareGeneralizationNonAblelian} or general approach to (dis)continuous cases \cite{HarrisonPoincareLemma, HarrisonOperatorCalculus}. In this paper, we will focus on the 'operator' approach to this lemma as well as on extension to complex manifolds.

As an introduction, we review some basic facts about the Poincar\'{e} lemma in order to fix notation. There are various formulations of the lemma and the most general one is the following well-known form 
\begin{Theorem} (Corollary 4.1.1 of \cite{DifferentialFormsInAlgebraicTopology}) (The Poincar\'{e} lemma)\\
 $$H^{*}(\mathbb{R}^n) = H^{*}(point)=\left\{ \begin{array}{l}
                                       \mathbb{R}, \quad (n=0) \\
                                       0 \quad (n>0)
                                      \end{array}\right.$$
\end{Theorem}
It can be formulated in another way by introducing (open) star-shaped region $U$ of $\mathbb{R}^{n}$ with respect to $x_{0}\in U$. It is an open region ($dim(U)=n$) where any other point $x\in U$ can be connected with $x_{0}$ using the line segment that lies entirely inside $U$. For a smooth manifold $M$ without boundary and of dimension $n=dim(M)$, we define a star-shaped region $U$ as (see, e.g., \cite{SmoothManifoldsLee}) the region diffeomorphic to an open ball in $\mathbb{R}^{n}$. For $M$ having boundary the star-shaped region can be also diffeomorphic to an open half-ball. Moreover, locally for a smooth manifold each point has a neighbourhood that is star-shaped.

Then we have
\begin{Theorem}(e.g. Theorem 11.49 of \cite{SmoothManifoldsLee}) (The Poincar\'{e} lemma) \\
 If $U$ is a star-shaped open subset of $\mathbb{R}^n$, then every closed covector field on $U$ is exact.
\end{Theorem}
As an existential statement, it is not useful in computations. The other approach relies on homotopy operator which is a local notion. To begin with, introduce the operator
\begin{equation}
 G\omega := \int_{0}^{1} (\partial_{t} \lrcorner \omega) dt,
\end{equation}
for $\omega \in \Omega(M \times \mathbb{R})$, where $\Omega(M \times \mathbb{R})$ is the module of forms, and where $M$ is a smooth manifold with or without boundary.
Next, choose a homotopy $F:[0,1]\times M \rightarrow M$ between $f$ and $g$, that is, $F(0,.)=f(.)$ and $F(1,.)=g(.)$. Using the homotopy we can define the operator
\begin{equation}
 \tilde{H}\omega = G\circ F^{*} (\omega),
 \label{Eg.HomotopyGeneral}
\end{equation}
for $\omega \in \Omega(M)$.
This operator has important property, which can be introduced using the Homotopy Invariance Formula, namely,
\begin{Theorem}(Paragraph 29 of \cite{ManifoldsTu})(Homotopy Invariance Formula for the de Rham complex)
\begin{equation}
dG + Gd = i_{1}^{*}-i_{0}^{*},
\label{Eq.HomotopyInvarianceFormula}
\end{equation}
where $i_{t}(x)=(t,x)$ for $t \in \mathbb{R}$ and $x \in M$.
\end{Theorem}
Using this formula we have the well known identity
\begin{equation}
(\tilde{H}d+d\tilde{H})\omega = GdF^{*}\omega+dGF^{*}\omega=i_{1}^{*}F^{*}\omega-i_{0}^{*}F^{*}\omega=g^{*}\omega-f^{*}\omega.
\label{Eq.GeneralHomotopyInvarianceFormula}
\end{equation}

The crucial observation was made by D.G.B. Edelen (see \cite{EdelenExteriorCalculus, EdelenIsovectorMethods}), that after a special choice of the homotopy, the definition (\ref{Eg.HomotopyGeneral}) has a particularly simple form with deep implications.  In order to derive Edelen's version of homotopy operator $H$ one have to choose the homotopy between the identity ($g(x)=x$) and the constant map ($f(x)=x_{0}$) for some fixed point $x_{0}\in U \subset M$. To provide the correct definition it is assumed that $U$ is a star-shaped region with respect to $x_{0}$. For such homotopy Edelen rewrote (\ref{Eq.GeneralHomotopyInvarianceFormula}) \cite{EdelenExteriorCalculus} as
\begin{Definition} (Edelen's homotopy operator)
 \begin{equation}
  H\omega := \int_{0}^{1} \mathcal{K} \lrcorner \omega_{F(t,x)} t^{k-1} dt,
  \label{Eq.EdelenHomotopyOperator}
 \end{equation}
for a $k$-form $\omega \in \Omega^{k}(U)$, $\mathcal{K}:=(x-x_{0})^{i}\partial_{i}$, $k=deg(\omega)$, and $F(t,x)=x_{0}+t(x-x_{0})$ is a homotopy between the constant map $x\rightarrow x_{0}$ and the identity map $I:x\rightarrow x$. The form $\omega$ under the integral is evaluated at the point $F(t,x)$.
\end{Definition}
The form of the operator is a special case of $\tilde{H}$ for the homotopy $F(t,x)=x_{0}+t(x-x_{0})$ and its explicit derivation is simple application of the pullback along $F$. $H$ has various properties described by Theorem 5-3.1 of \cite{EdelenExteriorCalculus}, from which the most important in later use is its nilpotency, $H^{2}=0$, which results from the double application of the insertion of $\mathcal{K}$ under the integral of (\ref{Eq.EdelenHomotopyOperator}).

The operator $H$ has its own Homotopy Invariance Formula, which can be written in a more compact form than in \cite{EdelenExteriorCalculus}, namely,
\begin{Theorem}(Homotopy Invariance Formula for $H$ operator)
\begin{equation}
 dH+Hd=I^{*} - s_{x_{0}}^{*},
 \label{Eq.EdelenHomotopyInvarianceFormula}
\end{equation}
where $s_{x_{0}}(x)=x_{0}$ is a constant map and $I$ is the identity.
\end{Theorem}
This formula was provided in Theorem 5-3.1 of \cite{EdelenExteriorCalculus} as a piecewise definition
\begin{equation}
 \left\{\begin{array}{c}
         Hd+dH=I, \quad \mathrm{on} \quad \Omega^k, k>0, \\
         (Hdf)(x) = f(x)-f(x_{0}) \quad \mathrm{for} \quad f \in \Omega^{0}.
        \end{array}
\right.
\end{equation}
It results from the fact that the pullback along the constant function $s_{x_{0}}^{*}\omega =0$ for $deg(\omega)>0$, and from the fact that $\mathcal{K} \lrcorner f =0$.

One can also note that (\ref{Eq.EdelenHomotopyInvarianceFormula}) is correct for any, not necessarily linear, homotopy $F$ between the identity and the constant map, however, in such a case, the explicit formula (\ref{Eq.EdelenHomotopyOperator}) is not valid.

If a form $\omega$ fulfils $d\omega = 0$ then it is called \textit{closed}. It is a well-known fact that in the star-shaped region $U$ (which we will assume hereafter), by the Poincar\'{e} lemma, it is also \textit{exact}, which means that there is a form $\alpha$ of degree $deg(\alpha)=deg(\omega)-1$ such that $\omega=d\alpha$. The exact (and hence closed) forms form a subspace $\mathcal{E}(U)$ of $\Omega(U)$. 

The following Lemmas will be useful in formulating operator calculus for $d$ and $H$:
\begin{Lemma}\label{Lemma.Antiexact1}(Lemma 5-4.1 of \cite{EdelenExteriorCalculus}) \\
 The operator $dH$ maps $\mathcal{E}^{k}$ onto $\mathcal{E}^{k}$ and $\Omega$ onto $\mathcal{E}$.
\end{Lemma} 
\begin{Lemma}\label{Lemma.Antiexact2}(Lemma 5-4.2 of \cite{EdelenExteriorCalculus}) \\
 The operator $d$ is the inverse of the operator $H$ when the domain of $H$ is restricted to $\mathcal{E}^{k}$.
\end{Lemma}
In addition, $\mathcal{E}^{0}(U)$ - the set of exact functions over $U$ is empty.

The less-known fact \cite{EdelenExteriorCalculus} is that, the $H$ singles out the so-called, \textit{antiexact} forms, that are the image of the complementary projection operator $Hd = I^{*}-dH-s_{x_{0}}^{*}$. This means that for an antiexact form $\omega$ there is an exact form $\alpha=d\beta$ such that $\omega = H\alpha$. The antiexact forms compose into the submodule $\mathcal{A}(U)$ of $\Omega(U)$ which is characterized by
\begin{Lemma}\label{Lemma.Antiexact3}(Lemma 5-5.1 of \cite{EdelenExteriorCalculus}) \\
\begin{equation}
\mathcal{A}^{k}=\{ \alpha \in \Omega^{k}, \mathcal{K}\lrcorner\alpha =0, \quad \alpha|_{x_{0}}=0, k>0\}. 
\end{equation}
\end{Lemma}
In addition, $\mathcal{A}^{n}(U)$ for $n=dim(U)$ is the empty set. Antiexactness is a local notion on star-shaped regions.

In this paper, we start from these simple properties and build on them additional abstract structures. The first aim is to formulate operator calculus in terms of $d$ and $H$. It is suitable to use the work of R. Bittner \cite{R.Bittner} who generalized differentiation and integration operations. This will allow us to formulate abstract differential equations and eigenvalue problems. We will show an example of such equations, which behaves similarly to the fermionic quantum harmonic oscillator (see, e.g., Chapter 5 of \cite{DasQFT}) used in quantum mechanics.

The second aim is to show how the Dolbeault bicomplex on complex manifolds interplay with the complex-valued version of $H$, which will be defined.

In summary, our aims are as follows:
\begin{itemize}
 \item {Make a more detailed characterization of exact and antiexact complexes.}
 \item {Construct operator calculus, where the exterior derivative plays a role of a 'derivative' and the homotopy operator is an 'integral'. This allows us to construct and solve abstract differential equations in these terms.}
 \item {Construct the homotopy operator for complex manifolds and describe its action on the Dolbeault complex.}
\end{itemize}

The paper is organized as follows. In the next section, a detailed description of $\Omega(U)$ decomposition into the exact vector space and the antiexact module will be given. Then in the following section, we will present the connection of these formulas with operator calculus and the fermionic quantum harmonic oscillator. Next, we develop the theory of the homotopy operator for complex manifolds.

\section{Homotopy operator and (anti)exact forms}

In \cite{EdelenExteriorCalculus} the decomposition $\Omega=\mathcal{E} \oplus \mathcal{A}$ on some star-shaped region $U$ of a smooth manifold $M$ is stated, however, using Lemmas \ref{Lemma.Antiexact1}, \ref{Lemma.Antiexact2}, \ref{Lemma.Antiexact3} and above properties, we have a finer characterization of this decomposition

\begin{Corollary} \label{Corollary.HomotopyComplex} {\ } \\ 
$\Omega(U)$ for $n=dim(U)>0$ is decomposed into the direct sum of exact and antiexact parts with respect to grading in the way presented in Fig. \ref{Fig.DecompositionFull}. At each degree $k \geq 0$ there is  $\Omega^{k} = \mathcal{E}^k \oplus \mathcal{A}^k$. The relations $d^2=0 = H^2$ when moving along the arrows are also visible.
 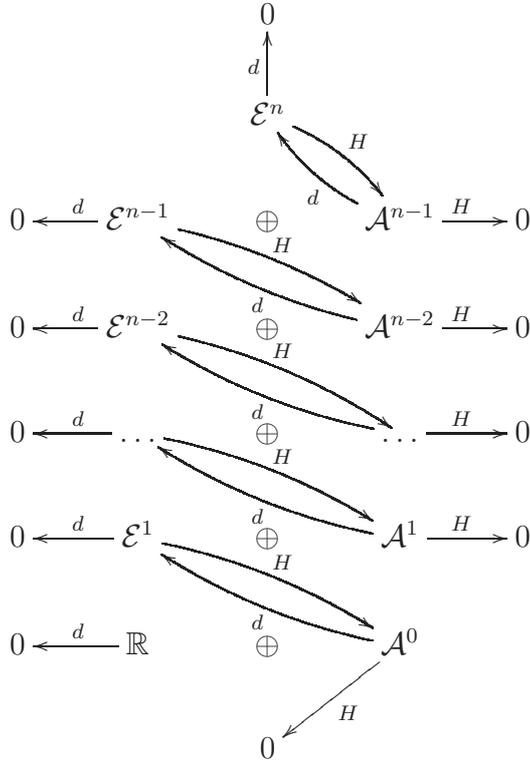
\begin{figure}
\centering
\xymatrix{ & &  0  & \\
       &     & \mathcal{E}^{n} \ar[u]^d \ar@/^/[dr]^{H} \\
0 & \ar[l]_d \mathcal{E}^{n-1} \ar@/^/[drr]^{H}  & \oplus & \mathcal{A}^{n-1} \ar[r]^H \ar@/^/[ul]^d  & 0 \\
0 & \ar[l]_d \mathcal{E}^{n-2} \ar@/^/[drr]^{H}  & \oplus & \mathcal{A}^{n-2} \ar[r]^H \ar@/^/[ull]^d  & 0 \\
0  & \ar[l]_d \ldots   \ar@/^/[drr]^{H}        &  \oplus  & \ldots \ar[r]^H \ar@/^/[ull]^d & 0  \\
0 & \ar[l]_d \mathcal{E}^{1} \ar@/^/[drr]^{H}  & \oplus & \mathcal{A}^{1} \ar[r]^H \ar@/^/[ull]^d  & 0 \\
0 & \ar[l]_d \mathbb{R} & \oplus & \mathcal{A}^{0} \ar@/^/[ull]^d  \ar[dl]^H   &  \\
     &  &   0   &
}
\caption{Decomposition of $\Omega$ into exact and antiexact subspaces with respect to the degree. Here $\mathbb{R}$ is treated as a space of constant functions.}
\label{Fig.DecompositionFull}
\end{figure}
\end{Corollary}

It is easily visible from the figure that for fixed $0<k<n=dim(U)$ there is separate 'subdiagram' depicted in Fig. \ref{Fig.DecompositionPart}, which will be the starting point for construction of operator calculus in the next section.
\begin{figure}
\begin{center}
\xymatrix{ 0 & \ar[l]_d \mathcal{E}^{k} \ar@/^/[r]^{H} & \mathcal{A}^{k-1} \ar[r]^H \ar@/^/[l]^d  & 0 \\
}
\end{center}
\caption{Part of the decomposition from Fig.\ref{Fig.DecompositionFull} for $0<k<n$.}
\label{Fig.DecompositionPart}
\end{figure}

For $k=0$, the kernel $Ker(d)$ consists of constant functions with values in the field over which $\Omega$ is the vector space, i.e., $\mathbb{R}$. This field can also be treated as constant $0$-forms. They are closed although not exact, which is a peculiarity in this decomposition.

The general formula (\ref{Eq.EdelenHomotopyInvarianceFormula}) is the starting point for considering the operator algebra of $H$, $d$, $I$ and $s_{x_{0}}$ in terms of operator calculus of Bittner \cite{R.Bittner} which will be the subject of the next section.

\section{Bittner's operator calculus}

\subsection{General setup}

The Bittner's operator calculus \cite{R.Bittner} is a way to redefine derivative and integral in abstract  terms. It mimics the well-known formulas
\begin{equation}
 \frac{d}{dx} \int_{q}^{x} f(x')dx' = f(x), \quad \int_{q}^{x} \frac{df(x')}{dx'} dx' = f(x)-f(q),
 \label{Eq.OperatorCalculusStandard}
\end{equation}
for $f$ being e.g. $C^{1}$ function. It is defined as follows
\begin{Definition}\cite{R.Bittner}\label{Def.BittnersOperatorCalculus} \\
Consider two linear spaces $L^{0}$ and $L^{1}$ and define an abstract derivative as surjective mapping $S \in Hom(L^{1},L^{0})$. Elements of $Ker(S)$ are called constants of the derivative $S$. Define also $T_{q} \in Hom(L^{0},L^{1})$ for some constant $q \in Ker(S)$ such that
\begin{equation}
 ST_{q}=I, \quad T_{q}S = I - s_{q},
 \label{Eq.BittnerOperatorCalculus}
\end{equation}
where $s_{q}$ is the projection operator on $Ker(d)$ associated with $q$. $T_{q}$ is called an abstract integral.
\end{Definition}
For instance, in (\ref{Eq.OperatorCalculusStandard}) $s_{q}f = f(q)\in ker\left( \frac{d}{dx}\right)$ is understood as a constant function.

We want to underline that this is not the derivation occurring in Differential Geometry \cite{SmoothManifoldsLee, Nakahara} since, e.g., no Leibnitz rule is implemented in this definition. Instead, it is a set of three operators that fulfill the requirements (\ref{Eq.BittnerOperatorCalculus}) that generalize (\ref{Eq.OperatorCalculusStandard}) from standard Calculus. It would be better to call them Bittner's derivative, integral, and projection on boundary conditions to distinguish them, however, we will not adopt this convention.

Let us consider the diagram from Fig. \ref{Fig.DecompositionPart} for $k>1$. In this case $Hd=I$ on $\mathcal{A}^{k-1}$ and $dH=I$ on $\mathcal{E}^{k}$ and therefore there is no projection on boundary data. In this case, restricted $H$ and $d$ are inverses to each other.

For $k>0$ this can be also seen as a mapping between the spaces on Fig. \ref{Fig.DecompositionPart2}.
\begin{figure}
\centering
\xymatrix{ 
0 & \ar[l]_d \mathcal{E}^{k} \ar@/^/[drr]^{H}  \\
0 & \ar[l]_d \mathcal{E}^{k-1}  & \oplus & \mathcal{A}^{k-1} \ar[r]^H \ar@/^/[ull]^d  & 0 \\
}
\caption{Part of decomposition of $\Omega^{k-1}=\mathcal{A}^{k-1}\oplus \mathcal{E}^{k-1}$ for fixed $0<k<n$. Note that $ker(d)=\mathcal{E}^{k-1}$ and $im(d) = \mathcal{E}^{k}$.}
\label{Fig.DecompositionPart2}
\end{figure}
In this case at the level $k-1$ we have $\Omega^{k-1}=\mathcal{A}^{k-1}\oplus \mathcal{E}^{k-1}$. This is decomposition according to the action of $d$ since $Ker(d)=\mathcal{E}^{k-1}$ and $Im(d)=\mathcal{E}^{k}$. We therefore have $Hd=I^{*}$ and $dH=I^{*}$, which is the special case of the second formula of (\ref{Eq.BittnerOperatorCalculus}) where $s_{x_{0}}^{*}=0$.

For $k=1$ for the case in Fig. \ref{Fig.DecompositionPart2} the pullback along the constant function $s_{x_{0}}^{*}$ is indispensable, and therefore, the formula (\ref{Eq.EdelenHomotopyInvarianceFormula}) leads to
\begin{equation}
 dH = I, \quad Hd = I^{*} - s_{x_{0}}^{*}.
 \label{Eq.BittnerHomotopyK=0}
\end{equation}
In this case the resemblance to (\ref{Eq.BittnerOperatorCalculus}) is even closer with $s_{q} = s_{x_{0}}^{*}$. In case of $\Omega^{0}$ the 'constants' of $d$ are the constant functions with values in $\mathbb{R}$, and pullback projects on them any function from $\Omega^{0}$ since it is evaluation of a function at $x=x_{0}$ and $\mathcal{A}^{0}\in ker(s_{x_{0}}^{*})$.

Therefore we have,
\begin{Lemma} {\ } \\ 
 For $k>0$ the operator calculus of Definition \ref{Def.BittnersOperatorCalculus} with the abstract derivative $d$ and the abstract integral $H$ is realized on the spaces of Fig. \ref{Fig.DecompositionPart2} as 
 \begin{equation}
  Hd=I \quad \mathrm{and} \quad dH=I^{*}.
 \end{equation}
For $k=0$ the Bittner's operator calculus is realized by (\ref{Eq.BittnerHomotopyK=0}).
\end{Lemma}

This observation allows us to formulate abstract differential equations on $\Omega$ using $d$ and $H$ as a 'derivative' and an 'integral' respectively. However, these operations are nilpotent, which put additional constraints on this 'operator calculus' and suggest that they resemble the situation appearing in the fermionic harmonic oscillator. This idea will be followed in the next subsection.

The idea of applying a combination of the exterior derivative and the homotopy operator to variational differential equations was discussed in \cite{IwaniecLutoborski} and \cite{DingNolder}, and this presentation uses methods of functional analysis. However, we follow a different direction focusing on connecting Edelen's idea and abstract Bittner's calculus.

\subsection{Homotopical harmonic oscillator}
We first recall a basic structure of the fermionic quantum harmonic oscillator. Consider a two-dimensional Hilbert space over $\mathbb{C}$. Then the fermionic quantum harmonic oscillator is defined by the Hamilton operator \cite{DasQFT}
\begin{equation}
 \hat{H}=a^{\dagger}a-aa^{\dagger},
 \label{Eq.FerminicHarmonicOscillator}
\end{equation}
where creation $a^{\dagger}$ and annihilation $a$ operator fulfills the anticommutation rules ($\{A,B\}:=AB+BA$)
\begin{equation}
 \{a,a\}=0, \quad \{a^{\dagger}, a^{\dagger}\}=0, \quad \{a, a^{\dagger} \}= I.
 \label{Eq.anticommutationRules}
\end{equation} 
The standard representation is
\begin{equation}
 a = \left( \begin{array}{cc}
             0 & 1 \\
             0 & 0
            \end{array} \right), \quad
a^{\dagger}=\left( \begin{array}{cc}
             0 & 0 \\
             1 & 0
            \end{array} \right).
\end{equation}
In this representation $\hat{H}$ is diagonal with eigenvalues $\pm1$.

The algebra of $d$ and $H$ is the same\footnote{We neglect the fact that $a^{\dagger}$ is an operator which is not the adjoint of $a$. In the space where $d$ and $H$ acts, there is no inner product which can be used to form such adjoint.} as for $a$ and $a^{\dagger}$, namely,
\begin{equation}
 dd=0, \quad HH=0, \quad Hd+dH = I^{*}-s_{x_{0}}^{*},
\end{equation}
where the term $s_{x_{0}}^{*}$ is zero when $deg(\omega)=k>0$.
It is therefore natural, by analogy, to define the Hamiltonian operator for the 'homotopical' harmonic oscillator
\begin{equation}
 \bar{H}:=Hd-dH.
 \label{Eq.FermionicHarmonicOscilatorHomotopyOpeator}
\end{equation}

We can now solve the eigenvalue problem for (\ref{Eq.FermionicHarmonicOscilatorHomotopyOpeator}), namely,
\begin{equation}
 \bar{H}\omega = \lambda \omega,
\end{equation}
where $\lambda \in  \mathbb{R}$ and $\omega \in \Omega^{k}$.
We have to consider three cases:
\begin{itemize}
 \item{$0<k<n$: The equation (\ref{Eq.FermionicHarmonicOscilatorHomotopyOpeator}) is of the form
 \begin{equation}
  2Hd\omega = (\lambda +1)\omega,
 \end{equation}
 and we are left with two cases:
 \begin{itemize}
  \item {$\lambda = -1$: for which $Hd\omega = 0$ that is $\omega \in  Ker(Hd)$, which gives that $\omega \in \mathcal{E}^{k}$.}
  \item {$\lambda \neq -1$: since $Hd\omega \in \mathcal{A}^{k}$ so $\omega \in \mathcal{A}^{k}$. Therefore $Hd\omega=\omega$ and the equation (\ref{Eq.FermionicHarmonicOscilatorHomotopyOpeator}) is $2\omega=(\lambda+1)\omega$, which gives $\lambda=1$.}
 \end{itemize}
 } 
 \item{$k=0$: take $f\in \Omega^{0}$, then $Hf=0$. If $f \in \mathcal{E}^{0}=ker(d)$ is a constant function then the eigenvalue problem for (\ref{Eq.FermionicHarmonicOscilatorHomotopyOpeator}) has the trivial solution $f=0$. Therefore we assume that $f\in  \mathcal{A}^{0}$. Then $Hdf = f-f_{x_{0}}$ and the eigenvalue problem  is 
 \begin{equation}
  f-f_{x_{0}}=\lambda f \quad \Leftrightarrow \quad (1-\lambda)f=f_{x_{0}}.
 \end{equation}
 For $f \in  \mathcal{A}^{0}$ we have $f_{x_{0}}=0$ and there are two cases
 \begin{itemize}
  \item {$\lambda =1$: then $f$ is an arbitrary element of $\mathcal{A}^{0}$.}
  \item {$\lambda \neq 1$: then $f=0$.}
 \end{itemize}
 }
 \item{$k=n$: let $\mu \in \mathcal{E}^{n}$, then $d\mu =0$ and $dH\mu=\mu$. The eigenvalue problem for  (\ref{Eq.FermionicHarmonicOscilatorHomotopyOpeator}) has the form
 \begin{equation}
  -dH\mu=\lambda \mu \quad \Leftrightarrow \quad (\lambda+1)\mu=0,
 \end{equation}
  which gives two cases:
  \begin{itemize}
   \item {$\lambda=-1$: then $\mu \in \mathcal{E}^{n}$ is arbitrary.}
   \item {$\lambda \neq -1$: then $\mu=0$.}
  \end{itemize}
 }
\end{itemize}

The above computation shows that the homotopical harmonic oscillator for $0<k<n$ only picks exact ($\lambda = -1$) or antiexact ($\lambda=1$) form and does not impose additional conditions on their functional form. For $k=0$ and $k=n$ there is only antiexact or exact solution respectively. This is a result of the fact that the tower from Fig. \ref{Fig.DecompositionFull} has a deficiency at the top and the bottom. It is an additional obstacle in making the analogy to the quantum mechanical system (\ref{Eq.FerminicHarmonicOscillator}) and (\ref{Eq.FermionicHarmonicOscilatorHomotopyOpeator}).

As in quantum mechanics \cite{DasQFT} there is also a top-down method for base generation where $H$ rises eigenstate and $d$ lowers eigenstate in the following sense:
\begin{itemize}
 \item {Let $\omega \in \mathcal{E}^{k}$, $k>0$. Then (since $d\omega=0$) locally $\omega=d\mu$ for $\mu \in \mathcal{A}^{k-1}$. Then $\hat{H} \omega= -dHd\omega= -d\mu = -\omega$, where the property $dHd=d$ from \cite{EdelenExteriorCalculus} was used. Therefore such $\omega$ is  $\lambda=-1$ eigenvector. It originates from $\mu$ which is $\lambda=1$ eigenvector.}
 \item {Likewise, let $\omega \in \mathcal{A}^{k}$, $k<n$. Then $H\omega=0$ and therefore $\omega = H\mu$. Finally, $\hat{H}\omega = HdH\mu = H\mu=\omega$, where the property $HdH=H$ of was used. Therefore $\omega$ is an eigenvector to the eigenvalue $\lambda = 1$. It originates from $\mu$ for eigenvalue $-1$.}
\end{itemize}
These two cases completely describe the diagram from Fig. \ref{Fig.DecompositionPart} and show how starting from one eigenvalue obtain the remaining one.

The homotopical fermionic quantum harmonic oscillator is, in some sense, similar to the Laplace-Beltrami operator known from Riemannian geometry \cite{Nakahara}. However, in our case, there is no metric structure and, as a result, no Hodge star, and therefore no codifferential can be constructed in a natural way. The homotopy operator is treated as a (local) alternative to the adjoint operator to the exterior derivative - the role that is played by the codifferential on a Riemannian manifold.

\section{Homotopy operator for complex manifolds}
This section contains an extension of the above theory for complex manifolds. We use the fact that Edelen's homotopy operator (\ref{Eq.EdelenHomotopyOperator}) does not 'feel' the underlying field of numbers.

First, we summarize the facts about complex manifolds in order to fix notation. Complex manifold \cite{Nakahara} is a smooth even dimensional manifold $M$ with holomorphic structure (of transition maps between coordinate patches). Such manifold has a complex structure $J$ which eigenspaces define the split of tangent space $T_{p}M = T_{p}M^{+} \oplus T_{p}M^{-}$, where the $+$ denotes the space spanned by holomorphic vector fields with the base $\{\partial_{z^{\mu}}\}_{\mu=1}^{n}$ and the space $-$ is spanned by anti-holomorphic vector fields with the base $\{\partial_{\bar{z}^{\mu}}\}_{\mu=1}^{n}$, where $2n = dim(M)$. We have
\begin{equation}
 \partial_{z^{\mu}}:=\frac{1}{2}\left( \partial_{x^{\mu}}-i \partial_{y^{\nu}}\right), \quad  \partial_{\bar{z}^{\mu}}:=\frac{1}{2}\left( \partial_{x^{\mu}}+i \partial_{y^{\nu}}\right),
\end{equation}
where $\{z^{1},\ldots,z^{n},\bar{z}^{1},\ldots,\bar{z}^{n}\}$ and $\{x^{1},\ldots,x^{n},y^{1},\ldots,y^{n}\}$ are local complex and real coordinates related by the standard formula $z^{\mu}=x^{\mu}+iy^{\mu}$.

This induces similar structure on the cotangent space, where the dual base has the $n$-dimensional covector base $dz^{\mu}$ of bidegree $(1,0)$ and the covector base $d\bar{z}^{\mu}$ of bidegree $(0,1)$. This constitutes the base of $1$-forms $\Omega^{1}(M)=\Omega^{1,0}(M)\oplus \Omega^{0,1}(M)$. Using exterior product, higher bidegree spaces can be constructed.

The exterior derivative $d$ can be decomposed as $d=\partial + \bar{\partial}$, where the Dolbeault operators are defined as
\begin{equation}
\begin{array}{cc}
 \partial:\Omega^{p,q}(M)\rightarrow \Omega^{p+1,q}(M), & \partial := dz^{\mu} \wedge \frac{\partial}{\partial z^{\mu}}, \\
 \bar{\partial}:\Omega^{p,q}(M)\rightarrow \Omega^{p,q+1}(M), & \bar{\partial} := d\bar{z}^{\mu} \wedge \frac{\partial}{\partial \bar{z}^{\mu}}.
\end{array} 
\end{equation}
Since from $d^{2}=0$ it results that $\partial^{2}=0$, $\bar{\partial}^{2}=0$ and $\partial\bar{\partial}+ \bar{\partial}\partial=0$ therefore they define a double complex on $\Omega^{p,q}(M)$.

Selecting a star-shaped region $U \subset M$ we can inside define, by analogy to (\ref{Eq.EdelenHomotopyOperator}), the homotopy operator where now $\mathcal{K} := (x-x_{0})^{\mu}\partial_{x^{\mu}}+(y-y_{0})^{\mu}\partial_{y^{\mu}}$, and the homotopy is $F(t,x,y)^{\mu}:=(x_{0}^{\mu}+t(x-x_{0})^{\mu}, y_{0}^{\mu}+t(y-y_{0})^{\mu})$. It is however more instructive to reformulate $H$ in terms of $z^{\mu}$ and $\bar{z}^{\mu}$ variables. In this case
\begin{equation}
 \mathcal{K}=\mathcal{K}^{+}+\mathcal{K}^{-},
\end{equation}
where 
\begin{equation}
 \mathcal{K}^{+}=(z-z_{0})^{\mu}\partial_{z^{\mu}}, \quad \mathcal{K}^{-}=\overline{\mathcal{K}^{+}}.
\end{equation}
Then the homotopy is $F(t,z)=z_{0}+t(z-z_{0})$ and similar for its complex conjugate. In this setup we have
\begin{Proposition}{\ } \\ 
$H$ splits into
\begin{equation}
 H=H^{+}+H^{-},
\end{equation}
where
\begin{equation}
 H^{\pm}\omega = \int_{0}^{1}\mathcal{K}^{\pm}  \lrcorner \omega_{F(t,z)} t^{k-1} dt.
 \label{Eq.ComplexEdelenHOperators}
\end{equation}
These operators act as follows
\begin{equation}
 H^{+}:\Omega^{p,q}(U) \rightarrow \Omega^{p-1,q}(U), \quad H^{-}:\Omega^{p,q}(U) \rightarrow \Omega^{p,q-1}(U),
\end{equation}
which vanish when $p-1<0$ or $q-1<0$, respectively. 

Then similarly to $H$ we have obvious properties
\begin{equation}
 H^{+}H^{+}=0 = H^{-}H^{-},
\end{equation}
and 
 \begin{equation}
  H^{+}H^{-}+H^{-}H^{+}=0,
 \end{equation}
which result from $HH=0$. 
\end{Proposition}

As a conclusion from the above Theorem and Corollary \ref{Corollary.HomotopyComplex} we have
\begin{Corollary}{\ } \\ 
$H^{\pm}$ define a double complex dual to the Dolbeault complex on a start-shaped region of a complex manifold. The complex is visualized in Fig. \ref{Fig.DolbeaultComplex}.
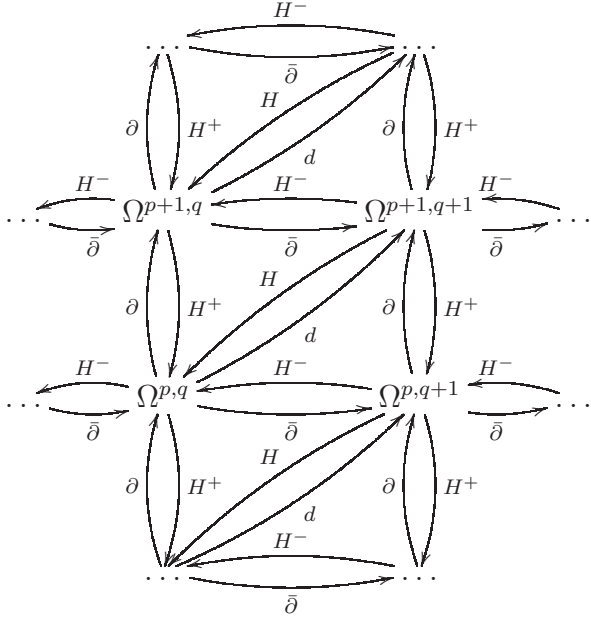
\begin{figure}
\centering
\xymatrix{ & \ldots \ar@/^/[dd]^{H^{+}} \ar@/_/[rr]_{\bar{\partial}} & & \ar@/_/[ll]_{H^{-}} \ldots \ar@/^/[dd]^{H^{+}} \ar@/_/[ddll]_{H}  & \\  & & & & &\\
\ldots \ar@/_/[r]_{\bar{\partial}}   &  \ar@/_/[l]_{H^{-}} \Omega^{p+1,q} \ar@/_/[rr]_{\bar{\partial}}  \ar@/^/[uu]^{\partial} \ar@/^/[dd]^{H^{+}} \ar@/_/[uurr]_{d}  & & \ar@/_/[ll]_{H^{-}} \Omega^{p+1,q+1} \ar@/_/[r]_{\bar{\partial}} \ar@/^/[uu]^{\partial} \ar@/^/[dd]^{H^{+}} \ar@/_/[ddll]_{H}  & \ar@/_/[l]_{H^{-}} \ldots  \\  & & & & &\\
 \ldots \ar@/_/[r]_{\bar{\partial}}   &  \ar@/_/[l]_{H^{-}} \Omega^{p,q}  \ar@/_/[rr]_{\bar{\partial}}  \ar@/^/[uu]^{\partial} \ar@/^/[dd]^{H^{+}} \ar@/_/[uurr]_{d} & & \ar@/_/[ll]_{H^{-}} \Omega^{p,q+1} \ar@/_/[r]_{\bar{\partial}} \ar@/^/[uu]^{\partial} \ar@/^/[dd]^{H^{+}} \ar@/_/[ddll]_{H}  & \ar@/_/[l]_{H^{-}} \ldots  \\  & & & & &\\
 &  \ldots \ar@/_/[rr]_{\bar{\partial}}  \ar@/^/[uu]^{\partial} \ar@/_/[uurr]_{d} &  & \ar@/_/[ll]_{H^{-}} \ldots  \ar@/^/[uu]^{\partial} &
}
\caption{Dolbeault complex and its homotopy dual.}
\label{Fig.DolbeaultComplex}
\end{figure}
\end{Corollary} 
 
In the complex case the formula (\ref{Eq.EdelenHomotopyInvarianceFormula}) becomes more elaborate
\begin{equation}
\begin{array}{c}
 Id-s_{(z_{0},\bar{z_{0}})}^{*} = (H^{+}+H^{-})(\partial+\bar{\partial})+(\partial+\bar{\partial})(H^{+}+H^{-})  \\
 =(H^{+}\partial + \partial H^{+}) + (H^{-}\bar{\partial} + \bar{\partial} H^{-}) + (H^{-}\partial + \partial H^{-}) + (H^{+}\bar{\partial} + \bar{\partial} H^{+}).
\end{array}
\label{Eq.CoplexEdelenHomotopyInvarianceFormula}
\end{equation}
The formula (\ref{Eq.CoplexEdelenHomotopyInvarianceFormula}) in general cannot be simplified to corresponding formulas for the pairs $(\partial, H^{+})$ and $(\bar{\partial}, H^{-})$ as it is presented in the following example. Consider a differential $(1,0)$ form $\omega = \bar{z}dz$. Nonzero elements of (\ref{Eq.CoplexEdelenHomotopyInvarianceFormula}) are
\begin{equation}
 \begin{array}{c}
  \partial H^{+} \omega = (\bar{z_{0}}+\frac{1}{2}(\bar{z}-\bar{z_{0}}))dz, \\
  \bar{\partial} H^{+} \omega = \frac{1}{2}(z-z_{0})d\bar{z}, \\
  H^{+}\bar{\partial} \omega = -\frac{1}{2}(z-z_{0})d\bar{z}, \\
  H^{-}\bar{\partial} \omega = \frac{1}{2}(\bar{z}-\bar{z_{0}})dz.
 \end{array}
\end{equation}
Summing these terms up we get $(Hd+dH)\omega = \bar{z} dz = I(\bar{z}dz) - s_{z_{0},\bar{z_{0}}}^{*}(\bar{z}dz)$ as required. Therefore all ingredients of (\ref{Eq.CoplexEdelenHomotopyInvarianceFormula}) must be taken into account in the general case. 

There are however two important cases when there is a split into $H^\pm$ subcomplexes.
\begin{Corollary} {\ } \\ 
There are two subcomplexes for $H^{+}$ and $H^{-}$, namely,
\begin{itemize}
 \item {$\bar{\partial}\omega = 0$ (holomorphic), $\omega \in \Omega^{p,0}, p\in \mathbb{N}$ - with no $d\bar{z}$ terms in the local representation, that is, $\omega = \omega(z)_{\mu_{1},\ldots,\mu_{p}}dz^{\mu_{1}}\wedge\ldots\wedge dz^{\mu_{p}}$. In this case $H^{-}\omega=0$ (anti-$\bar{\partial}$-exact), and $\bar{\partial}H^{+}\omega=0$. Then (\ref{Eq.CoplexEdelenHomotopyInvarianceFormula}) has the simple form
 \begin{equation}
  H^{+}\partial + \partial H^{+}=I - s_{z_{0}}^{*}.
 \end{equation}
 
 This defines the subcomplex $(\Omega^{p,0}, \partial, H^{+})$.
 }
 \item {$\partial\omega = 0$ (antiholomorphic), $\omega \in \Omega^{0,p}, p\in \mathbb{N}$ - with no $dz$ terms in the local representation, that is, $\omega = \omega(\bar{z})_{\mu_{1},\ldots,\mu_{p}}d\bar{z}^{\mu_{1}}\wedge\ldots\wedge d\bar{z}^{\mu_{p}}$. In this case $H^{+}\omega =0$ (anti-$\partial$-exact), and $\partial H^{-}\omega=0$. Then (\ref{Eq.CoplexEdelenHomotopyInvarianceFormula}) has the simple form 
 \begin{equation}
  H^{-}\bar{\partial} + \bar{\partial} H^{-}=I - s_{\bar{z_{0}}}^{*}.
 \end{equation}
 Likewise, this defines the subcomplex $(\Omega^{0,p}, \bar{\partial}, H^{-})$.
 }
\end{itemize}

Both of these subcomplexes lie on the boundary (left and bottom part) of Fig. \ref{Fig.DolbeaultComplex}.
 
\end{Corollary}

\section{Conclusions}
In this paper, the local results related to the homotopy operator from the Poincar\'{e} lemma was used to derive a special case of operator calculus that resembles structures occurring in quantum mechanics. Moreover, the analysis of dual Dolbeault bicomplex induced on a complex manifold by complex homotopy operator was provided. Two special subcomplexes were identified. These results organize and generalize the Poincar\'{e} lemma by building additional abstract structure on the top of this classical result.

\section*{Acknowledgments}
I would like to thank Josef \v{S}ilhan for long stimulating discussions during the preparation of the paper, and Luk\'{a}\v{s} Vok\v{r}\'{\i}nek for discussion about Homological algebra. I would also like to thank Jan Slov\'{a}k and Henrik Winther for useful suggestions. Last but not least, I thank anonymous Referee, whose vital and precise comments and suggestions improve this paper.

This research was supported by the GACR grant GA19-06357S and Masaryk University grant MUNI/A/0885/2019. I also thank the PHAROS COST Action (CA16214), and SyMat COST Action (CA18223) for partial support.

\appendix




\end{document}